\documentclass{article}
\usepackage{aism-2e}
\usepackage{amsmath,amssymb}
\newtheorem{theorem}{\sc Theorem}

\newtheoremrm{definition}{\sc Definition}[thrmctr]
\newtheoremrm{assumption}{\sc Assumption}[thrmctr]
\newtheoremrm{condition}{\sc Condition}[thrmctr]
\newtheoremrm{remark}{\it Remark\/}                  
\newtheoremrm{example}{\it Example\/}
\newtheoremrm{case}{\it Case\/}
\newtheoremrmno{proof}{\sc Proof}
\newtheoremrmno{procedure}{\sc Procedure}
\newtheoremrmno{note}{\it Note\/}
\newtheoremrmno{fact}{\it Fact\/}
\makeatletter
\@addtoreset{corollary}{section}
\@addtoreset{proposition}{section}
\@addtoreset{assumption}{section}
\@addtoreset{condition}{section}
\frompage{1}\topage{12}

 \newcommand{\cE}{\mathcal{E}}

 \newcommand{\bC}{\mathbb{C}}
 \newcommand{\bR}{\mathbb{R}}
\newcommand{\bN}{\mathbb{N}}

\newcommand{\erf}{\mathrm{erf}}

\begin{document}
\sloppy

\title{
On the Laplace transform of some quadratic forms and the exact distribution of the sample variance from a gamma or uniform parent  distribution}

\author{
T. Royen
}
\affiliation{
Fachhochschule Bingen, University of applied sciences,\\
 Berlinstrasse 109, D--55411 Bingen, Germany\\
 E--mail: royen@fh-bingen.de
}


\abstract{
From a suitable integral representation of the Laplace transform of a positive semi--definite quadratic form of independent real random variables with not necessarily identical densities a univariate integral representation is derived for the cumulative distribution function of the sample variance of i.i.d. random variables with a gamma density, supplementing former formulas of the author. Furthermore, from the above Laplace transform  Fourier series are obtained for the density and the distribution function of the sample variance of i.i.d. random variables with a uniform distribution. This distribution can be applied e.g. to a statistical test for a scale parameter. 
} 

\keywords {
Quadratic forms of non--normal random variables, Laplace transform of quadratic forms, Exact distribution of sample variances, Gamma distribution, Uniform distribution, Test of scale

\noindent\textsl{AMS 2000 subject classifications:}  62E15, 62H10}
\date{2.10.2007}
\maketitle

\section*{\large 1. \ Introduction }
\label{section:1}
Exact results for the distribution of the sample variance $s^2$ from non--normal samples are hardly found. Early investigations, concerning $s^2$ from a sample of gamma random variables can be traced at least to Craig (1929) and Pearson (1929), who used moment approximations, later investigated more thoroughly by Bowman and Shenton (1983). For general distributions various approximations were proposed by Box (1953), Roy and Tiku (1962), Tan and Wong (1977) and by Mudholkar and Trivedi (1981). The latter authors recommend transformations of the Wilson--Hilferty type. For a mixture of two normal distributions see also chapter 4.7 in Mathai and Provost (1992). This monograph treats quadratic forms, mainly from normal random variables. For the cdf of a positive semi--definite quadratic form of normal r.v. a simple integral representation on $[0,\pi]$ is given by formula (31) in Royen (2007a) as the univariate special case of some multivariate distributions. The interest in the distribution of $s^2$ or $1/s$ was also stimulated by the estimation of process capability indices, see e.g. Pearn, Kotz and Johnson (1992) and Pearn and Kotz (2006).

In this paper the cumulative distribution function (cdf) of $(n-1)s^2 = \sum_{j=1}^n (X_j - \overline{X})^2$ is derived for i.i.d. random variables $X_j$ with a gamma density $x^{\alpha-1} e^{-x}/\Gamma(\alpha)$ and for uniformly distributed $X_j$. 
The first distribution was already given by a double series in Royen (2007b), but now a univariate integral representation of the cdf is derived in section~3 avoiding the recursive computation of the coefficients of the above series. The cited paper contains also the exact distribution of $\Big( \sum^n_{j=1} X_j^2\Big)^{1/2}\Big/\sum^n_{j=1} X_j$ including Greenwood's statistic if $\alpha = 1$.
 The distribution of $(n-1)s^2$ from a sample of uniformly distributed random variables is obtained by a Fourier series in section~4 and applied to a statistical test for a scale parameter. 

Both these distributions are derived from a suitable integral representation of the Laplace transform (LT) of an "$m$--factorial"positive semi--definite quadratic form
\begin{equation}
\label{eq:1}
Q = X' (D-CC')X
\end{equation}
of independent real random variables $X' = (X_1, \ldots, X_n)$ with not necessarily identical densities, with positive definite $D = \diag (d_1, \ldots,d_n)$ and an $n \times m$--matrix $C = (c_{j\mu})$ of rank $m$ with real or pure imaginary columns $c_\mu$.

In particular, with i.i.d. random variables and $m=1$ the form $Q$ is reduced to $(n-1)s^2$ if $D$ is the identity matrix and $c_j = 1/\sqrt n$ for all $j$. The LT of $Q$ is given in the following section.

Formulas from the handbook of mathematical functions by Abramowitz and Stegun are cited by A.S. and their number. The standard normal density is denoted by $\varphi$ and $\cE$ means "expectation".

\section*{\large 2. \ The Laplace transform of an $m$--factorial positive semi--definite quadratic form }
\label{section:2}

The LT or Fourier transform (FT) of a sum of independent random variables can be written as a product of univariate integral transforms. To preserve this advantage for a quadratic form of the type in (\ref{eq:1}) as far as possible its LT is represented by an $m$--variate mixture over a product of univariate integral transforms. This is accomplished essentially by the equation
\begin{equation}
\label{eq:2}
\int_{-\infty}^\infty \exp(yz) \varphi(y)dy = \exp(z^2/2), z\in \bC .
\end{equation}
Inspite of its simplicity the result is stated formally as a theorem.\\[1ex]

\begin{theorem}
\label{theorem:1}
Let $Q = X' (D-CC')X$ be a positive semidefinite quadratic form of the type described in (\ref{eq:1}). Then the LT of $Q$ is given by 
\begin{equation}
\label{eq:3}
\cE (\exp(-tQ)) = \int_{\bR^m} \Big( \prod_{j=1}^n \cE\big(\exp\big(-td_j X_j^2 \pm \sqrt{2t} X_j \sum_{\mu=1}^m c_{j\mu} y_\mu\big)\big)\Big) \cdot \prod_{\mu=1}^m \varphi(y_\mu) dy_\mu 
\end{equation}
with the standard normal density $\varphi, t \geq 0$, all $d_j > 0$.
\end{theorem}

\paragraph{\sc Remark: } Replacing $t$ by $-it$ in (\ref{eq:3}) the characteristic function of $Q$ is obtained if $Q$ is definite or indefinite and absolutely bounded.

With i.i.d. random variables, $m=1$, $d_j\equiv 1$, $c_j \equiv 1\sqrt n$ we get
\begin{eqnarray}
\label{eq:4}
\cE\left( \exp(-tQ)\right) &=& \cE \left( \exp(-t(n-1)s^2\right)\nonumber\\
&=& \int^{\infty}_{-\infty}	\big(\cE\big( \exp(-tX^2 \pm \sqrt{2t/n} Xy )\big)\big)^n \varphi(y) dy\\
&=& \frac 1{\sqrt{2\pi}} \int^{\infty}_{-\infty} \big(\cE\big( \exp(-(\sqrt{t} X \ {\overline{+}}\
y/\sqrt{2n} )^2)\big)\big)^n dy .\nonumber
\end{eqnarray}

\paragraph{\sc Proof of theorem 1 :} Let be $m =1$:
\begin{eqnarray*}
	\cE\Big( \exp(-tQ)\Big) &=& \cE \Big( \exp \big(-t \sum^{n}_{j=1} d_j X^2_j + t \big( \sum^{n}_{j=1} c_j X_j\big)^2\big)\Big)\\
&=& 	\cE\Big( \exp(-t \sum^{n}_{j=1} d_j X^2_j\big) \int^{\infty}_{-\infty} \exp \big( \pm y \sqrt{2t} \sum^{n}_{j=1} c_j X_j\big) \varphi (y) dy\Big), 
\end{eqnarray*}
applying (\ref{eq:2}) with $z = \pm \sqrt{2t} \sum_j c_j X_j$.

For real $c_j$ the integrand is only positive. Thus, equation (\ref{eq:3}) is obtained by changing the order of integration. For imaginary values $c_j$ the integrand is majorized by $\exp(-tX'DX)\varphi(y)$ and equation (\ref{eq:3}) follows in the same way. The generalization to $m > 1$ is obvious and needs no further considerations. 

\section*{\large 3. \ The exact distribution of the sample variance from a gamma parent \mbox{distribution} }
\label{section:3}

We need a moment generating function $M_\alpha$ and a "kernel" $K_\beta$:
 \begin{eqnarray}
\label{eq:5}
M_\alpha(z) &:=& (\Gamma(\alpha))^{-1} \int^{\infty}_{0} \exp(-x^2) x^{\alpha-1} \exp(xz)dx\\
&=& (2\Gamma(\alpha))^{-1} \sum^{\infty}_{k=0} \Gamma ((\alpha + k)/2) z^k/k!\nonumber\\
&=&(2\Gamma(\alpha))^{-1} \Big( \Gamma \big( \frac \alpha 2 \big)\ _1F_1 \big( \frac \alpha 2, \frac 12, \frac{z^2}4\big) 
+ \Gamma \big( \frac{\alpha +1} 2 \big) z \ _1F_1 \big( \frac{\alpha +1}2, \frac 32, \frac{z^2}4\big)\Big),\nonumber\\
&& \qquad \alpha > 0, z \in \bC.\nonumber
\end{eqnarray}
If $\alpha \in \bN$ then 
\begin{equation}
\label{eq:6}
M_\alpha(z) = (2\Gamma(\alpha))^{-1} \sqrt \pi \big(\frac d{dz}\big)^{\alpha-1} \left( \exp(z^2/4) \ \mathrm{erfc} (-z/2)\right).
\end{equation}
\begin{eqnarray}
\label{eq:7}
K_\beta(r,z) &:=& \sum^{\infty}_{k=0} \big( \Gamma(\beta + k/2)  \big)^{-1} He_k(z) (-r)^k /k!\\
&=& \sum^{\infty}_{k=0} \big( \Gamma(\beta + k/2)  \big)^{^-1} \, _0F_1 (\beta + k/2, -r^2/2) (-rz)^k/k!,\nonumber\\
&& \qquad \beta > 0, \ r \geq 0, \ z\in \bC, \nonumber
\end{eqnarray}
where $He_k(z) = (-1)^k \varphi^{(k)} (z) / \varphi(z) = 2^{-k/2} H_k (z/ \sqrt 2)$ with the Hermite polynomials $H_k$.

The equivalence of the two series in (\ref{eq:7}) follows by inversion of the LT
\begin{eqnarray}
\label{eq:8}
\lefteqn{t^{-\beta}\varphi(z+t^{-1/2})/\varphi(z) = t^{-\beta} \exp (-z/\sqrt t) \exp(-1/(2t))} \\
&=& \exp (-1/(2t)) \sum^{\infty}_{k=0} t^{-\beta -k/2} (-z)^k /k!
= \sum^{\infty}_{k=0} (-1)^k He_k(z) t^{-\beta -k/2}/k!\nonumber
\end{eqnarray}
with $t > 0$ and any fixed $z$. The inversion provides the two corresponding series for $x^{\beta -1} K_\beta (\sqrt x, z)$.

From (\ref{eq:7}) we obtain also an integral representation of $K_\beta$ by well known functions. Let be
\begin{eqnarray*}
w_\beta(y) &=& \sum^\infty_{k=0} (\Gamma(\beta + k/2))^{-1} y^k\\
&=& (\Gamma(\beta))^{-1} \  _1F_1(1,\beta,y^2) + (\Gamma(\beta + 1/2))^{-1} y\ _1F_1(1,\beta + 1/2,y^2), \ y \in \bC.
\end{eqnarray*}
After Kummer's transformation (A.S.13.1.27) we find with (A.S.13.6.10) and \mbox{$G_\beta(z) = (\Gamma(\beta))^{-1} \int_0^z x^{\beta-1} e^{-x} dx$}
\[
w_\beta(y) = \exp (y^2) y^{-2(\beta-1)} (G_{\beta-1}(y^2) + G_{\beta-1/2}(y^2)), \ \beta > 1.
\]
In particular with $\beta = 1+m\in \bN$ we have 
\begin{eqnarray*}
w_\beta(y) &=& y^{-2m}\Big( \exp (y^2) \ \mathrm{erfc}(-y) - \sum_{j=0}^{2m-1} y^j \big/ \Gamma(1+j/2)\Big),\\
w_{1/2}(y) &=& y \exp(y^2) \ \mathrm{erfc}(-y) + \pi^{-1/2},\\
w_{1/2+m}(y) &=& y^{-1} (w_m(y) - 1/(m-1)!).
\end{eqnarray*}
Thus, with any $\rho >0$, we obtain 
\[
K_\beta(r,z)\varphi(z) = \frac 1{2\pi} \int^\pi_{-\pi} \varphi(z + \rho e^{i\psi}) w_\beta(\rho^{-1} re^{-i\psi})d\psi.
\]
$K_\beta$ will be used only with real $z$. In this case the real part of the integrand is integrated on $[0,\pi]$.

Furthermore, from the $2^{nd}$ series in (\ref{eq:7}) and the integral representation of the Bessel functions $J_{\beta-1+k/2}$ it follows
\begin{eqnarray*}
\lefteqn{K_\beta(r,z) = (2/\sqrt\pi)\int_0^{\pi/2} \Big(\, _0F_2\big(\frac 12, \beta-\frac 12, \frac 14 \ r^2 z^2 \sin^2 \vartheta\big)\big/ \Gamma(\beta-\frac 12 )\ - }\\
&& (r z \sin \vartheta) \ _0F_2\big( \frac 32, \beta, \frac 14 \ r^2 z^2 \sin^2 \vartheta\big)\big/\Gamma(\beta)\Big)(\sin^2 \vartheta)^{\beta-1}\cos (r \sqrt 2 \cos \vartheta) d\vartheta, \ \beta > \frac 12 \ .
\end{eqnarray*}

For the proof of the following theorem 2 we need some knowledge concerning the behaviour of the functions $M_\alpha$ and $K_\beta$ for large $|z|$, $z$ real. From the integral in (\ref{eq:5}) and the Tauber theorem for Laplace transforms of measures (e.g. Feller, chapter 13) it follows
\begin{equation}
\label{eq:9}
\lim_{z\to -\infty} |z|^\alpha M_\alpha (z) = 1.
\end{equation}
With (A.S.9.1.62) and (A.S.9.1.69) we have $| _0F_1 (\beta + k/2, - r^2/2)| \leq 1$. Then, from $\lim_{k\to \infty} \varepsilon^{-k} (\Gamma (\beta + k/2))^{-1} =0$ for all $\varepsilon >0$, we obtain with the second series in (\ref{eq:7})
\begin{equation}
\label{eq:10}
K_\beta (r,z) = o(e^{\varepsilon r |z|}), \ \varepsilon > 0, \ z\to -\infty.
\end{equation}
With (A.S.22.14.17) we get $|He_k(z)|/k! < 1.1  e^{z^2/4}/\sqrt{k!}$ and consequently
\begin{equation}
\label{eq:11}
|K_\beta(r,z)|< 1.1 e^{z^2/4} \sum^{\infty}_{k=0}\big(\Gamma(\beta + k/2)\sqrt{k!}\big)^{-1} r^k, \ z \in \bR.
\end{equation}
From the last representation in (\ref{eq:5}) and (A.S.13.1.14) it follows the asymptotic relation
\begin{equation}
\label{eq:12}
M_\alpha (z) \simeq (\sqrt \pi/\Gamma(\alpha)) \exp (z^2/4)(z/2)^{\alpha-1}, \ z \to \infty.
\end{equation}
Finally, replacing $t$ by $-it$ in (\ref{eq:8}) we find with $z > 0$ and $\beta >1$ by the Fourier inversion formula
\begin{eqnarray}
\label{eq:13}
|x^{\beta -1} K_\beta (\sqrt x,z)| &\leq& \pi^{-1} \int^{\infty}_{0} t^{-\beta} \exp (-z /\sqrt {2t}) dt\\
&=& (2^\beta/\pi) \Gamma (2 \beta -2) /z^{2\beta -2}.\nonumber
\end{eqnarray}\\[1ex]

\noindent \begin{theorem}
\label{theorem:2}
Let $X_1,\ldots,X_n$ be i.i.d. random variables with a gamma density $x^{\alpha -1} \exp (-x) /\Gamma(\alpha)$. With the functions $M_\alpha$ from (\ref{eq:5}) and $K_\beta$ from (\ref{eq:7}) the cdf of $Q = \sum^{n}_{i=1} (X_i - \overline{X})^2$, $n \geq 2$, is given by
\begin{equation}
\label{eq:14}
\Pr\{Q \leq x\} = x^{\alpha n/2} \int^{\infty}_{-\infty} K_{\alpha n/2+1} (\sqrt{nx/2}, z) (M_\alpha (\sqrt{2/n}\  z))^n \varphi (z) dz.
\end{equation}
\end{theorem}

\begin{proof}
Due to the relations (\ref{eq:9})$\ldots$(\ref{eq:13}) with $\beta =\alpha n/2 + 1$ the integral in (\ref{eq:14}) is absolutely convergent for all $\alpha > 0$ and $x \geq 0$. In particular, the integrand is at least $O(z^{-n})$ for $z \to \infty$. Furthermore, the LT  
\[
\int^\infty _0 e^{-tx} x^{\alpha n/2} \big( \int^\infty_{-\infty} \big|K_{\alpha n/2+1} (\sqrt{nx/2}, z)\big| (M_\alpha(\sqrt{2/n} \ z))^n \varphi (z)dz\big)dx 
\]
 exists for all $t > 0$.  Thus, changing the order of integration, the LT of
\[
x^{\alpha n/2} \int^\infty_{-\infty} K_{\alpha n/2+1} (\sqrt{nx/2},z) \big(M_\alpha (\sqrt{2/n}\ z)\big)^n \varphi(z)dz
\]
is given by
\begin{equation}
\label{eq:15}
t^{-(\alpha n/2+1)}\int^\infty_{-\infty} \Big( \varphi(z + \sqrt{n/(2t)})/\varphi (z)\Big)\Big(M_\alpha\big( \sqrt{2/n}\ z\big)\Big)^n \varphi(z) dz.
\end{equation}
On the other hand, applying equation (\ref{eq:4}) to a sample of i.i.d. random variables with density $x^{\alpha -1} \exp (-x) / \Gamma (\alpha)$, we find  the LT of the density of $Q$.
\begin{eqnarray*}
\lefteqn{\cE (e^{-tQ}) = (\Gamma(\alpha))^{-n} \int^\infty_{-\infty} \Big( \int^\infty_0 \exp(-tx^2 + \sqrt{2t/n}\ xy-x)x^{\alpha-1} dx\Big)^n \varphi(y)dy}\\
&=& (\Gamma(\alpha))^{-n} t^{-\alpha n/2}  \int^\infty_{-\infty} \Big( \int^\infty_0 \exp(-x^2 + (\sqrt{2/n}\ y - 1/\sqrt t)x)x^{\alpha-1}dx\Big)^n \varphi(y) dy\\
&=& (\Gamma(\alpha))^{-n} t^{-\alpha n/2}  \int^\infty_{-\infty} \Big( \int^\infty_0 \exp(-x^2 + \sqrt{2/n}\ xz) x^{\alpha-1}dx\Big)^n \varphi\big(z+ \sqrt{n/(2t)}\big) dz
\end{eqnarray*}
with the substitution $z = y - \sqrt{n/(2t)}$. Therefore $t^{-1} \cE(e^{-tQ})$ coincides with (\ref{eq:15}) which completes the proof.
\end{proof}

In the same way the more general formula
\[
\Pr \{Q \leq x\} = x^{\alpha/2} \int^\infty_{-\infty} K_{\alpha/2+1} \big( \sqrt{nx/2}, z\big) \Big( \prod^n_{j=1} M_{\alpha_j} \big( \sqrt{2/n}\  z\big)\Big) \varphi(z)dz, \ \alpha = \sum^n_{j=1} \alpha_j
\]
is proved for independent $\Gamma(\alpha_j)$--distributed random variables $X_j$.

\section*{\large 4. \ The exact distribution of the sample variance from a uniform parent \mbox{distribution}}
\label{section:4}
Let $X_1,\ldots,X_n$ be i.i.d. random variables, uniformly distributed on the interval $[0,1]$. The maximum $Q_{\max}$ of $Q = \sum^n_{i=1}(X_i - \overline{X})^2$ is given by $n/4$ if $n$ is even and by $(n^2-1)/(4n)$ otherwise. These values are reached if $X_i(1-X_i) = 0$, $i=1,\ldots,n$, with exactly $[n/2]$ solutions $X_i = 1$ or $[n/2]$ solutions $X_i=0$. For a formal proof set $Y_i = 2X_i-1 = \sin \Phi_i$. The coordinates of the stationary points of $4Q = \sum^n_{i=1} \sin^2 \Phi_i - 1/n (\sum^n_{i=1}\sin \Phi_i)^2$ satisfy $Y_i = \pm 1$ or $Y_i = \overline{Y}$. It is easy to show that the maximum is only attained if $Y_i = \overline{Y}$ does not occur. Then the result follows by maximizing $Q = k(n-k)/n$, $k = 1,\ldots,n-1$. 

From (\ref{eq:4}) the LT of $Q$ follows by straightforward calculation. Hence, the Fourier transform $\widehat{f}$ of the density $f$ of $Q$ is given by
\begin{equation}
\label{eq:16}
\widehat{f}(t) = \sqrt n(\sqrt \pi/2)^{n-1}(i/t)^{n/2} \int^\infty_0 \Big(\erf\big(y+\frac 12\sqrt{-it}\big) - \erf\big(y-\frac 12 \sqrt{-it}\big)\Big)^n dy
\end{equation}
and the Fourier sinus coefficients $b_k$ of $f(x)$, $0 \leq x\leq Q_{\max}$, are 
\begin{eqnarray}
\label{eq:17}
b_k &=& (2/Q_{\max}) \Im m(\widehat{f}(t_k)),\nonumber\\[-2ex]
&&~  \\[-2ex]
t_k &=& \Bigg\{ \begin{matrix} 4k\pi/n & , \ n \ \mathrm{even}\\
                               4k\pi n/ (n^2-1) & , \ n \ \mathrm{odd} \end{matrix} \Bigg\}.\nonumber
\end{eqnarray}
However, for numerical reasons, any values of $\erf(z)$ with $|\Im m(z)| >|\Re e(z)|$ should be avoided within the integral $\widehat{f}(t)$. Therefore, the path of integration is shifted by $1/2\sqrt{-it}$. With Cauchy's integral theorem for holomorphic functions and 
\[
\lim_{y \to \infty} \erf (y \pm \sqrt{-it}) - \erf (y) = 0,
\]
uniformly on intervals $0 \leq t \leq T$, the integral 
\[
\int^\infty_0 \Big(\erf \big(y + \frac 12 \sqrt{-it}\big) - \erf \big(y - \frac 12 \sqrt{-it}\big)\Big)^n dy
\]
in (\ref{eq:16}) can be replaced by the identical value
\begin{eqnarray}
\label{eq:18}
\lefteqn{\int^\infty_0 \Big(\erf \big(y + \sqrt{-it}\big) - \erf(y)\Big)^n dy }\\
&& + \frac 12 \sqrt{-it} \int^1_0 \Big(\erf \big (\frac 12 \sqrt{-it}(u+1)\big)-\erf\big(\frac 12 \sqrt{-it} (u-1)\big)\Big)^n du.\nonumber
\end{eqnarray}
With the values $t_k$ from (\ref{eq:17}) this leads to\\[1ex]

\begin{theorem}
\label{theorem:3}
Let $X_1,\ldots,X_n$ be independent random variables, uniformly distributed on the interval $[0,1]$. The cdf of $Q = \sum^n_{i=1} (X_i - \overline{X})^2$ is given by

\begin{eqnarray}
\label{eq:19}
\lefteqn{F(x) = \Pr\{Q \leq x\}= \frac 12 \sqrt n (\sqrt \pi/2)^{n-3} \ \cdot \nonumber }\\
&& \Bigl( \int^\infty_0 \sum^\infty_{k=1} \Im m\Big( (i/t_k)^{n/2} \big(\erf(y+\sqrt{-it_k}) - \erf(y)\big)^n\Big)
(1-\cos(t_kx))/k\,dy \ + \\
&&  \frac 12\int^1_0  \sum^\infty_{k=1}\Im m\Big( (i/t_k)^{(n-1)/2} \big(\erf \big(\frac 12 \sqrt{-it_k}(u+1)\big) -\erf\big( \frac 12 \sqrt{-it_k}(u-1)\big)\big)^n\Big) \cdot\nonumber\\
&&(1-\cos (t_kx))/k\,du\Bigr), \nonumber
\end{eqnarray}
\[ 
0 \leq x < Q_{\max} = \Bigg\{ \begin{matrix} n/4 & \ , \ n \ \mathrm{even}\\ (n^2-1)/(4n) & \ , \ n \ \mathrm{odd} \end{matrix} \Bigg\} ,\quad  n \geq 3.
\]
The density $f(x)$ is obtained by the derivative under the integral at  least for $n \geq 4$.
\end{theorem}

\begin{proof}
The formal series 
\begin{equation}
\label{eq:20}
\frac 2\pi \sum^\infty_{k=1} \Im m(\widehat{f}(t_k)) (1- \cos (t_kx))/k
\end{equation}
for $F(x)$ follows by integration of the formal Fourier sinus series for the density $f$. Since $F(x)$ is geometrically the intersection volume of the unit $n$--cube with an $n$--cylinder of radius $r = \sqrt x$ whose axis passes  $(0,\ldots,0)$ and $(1,\ldots,1)$, we have $F(x) = O(x^{(n-1)/2})$ and $f(x) = O(x^{(n-3)/2})$ for $x \to 0_+$. Besides, $f(x)$ tends to zero for $x \to Q_{\max}$. Hence, $f$ is square integrable for all $n \geq 3$ and consequently the formula (\ref{eq:20}) provides an absolutely convergent series representation of $F(x)$. After the replacement described in (\ref{eq:18}) the resulting series differs from  (\ref{eq:19}) only by the change of summation and integration. Since $\erf(z)$ is absolutely bounded on $\{z \in \bC \big| |\Im m (z)| \leq |\Re e(z)|\}$ the change is justified for the $2^{nd}$ integral in (\ref{eq:19}). The first integral in (\ref{eq:19}) is divided into $\int^1_0 + \int^\infty_1$ and we get with $y \geq 1$ and $t = t_k > 0$ the estimation 
\begin{eqnarray*}
\lefteqn{\big|\big( \erf(y + \sqrt{-it}) - \erf(y)\big)/\sqrt{-it}\ \big| } \\
&=& \frac 2{\sqrt \pi} \big|(-it)^{-1/2} \int^{y + \sqrt{-it}}_y \exp(-z^2)dz\big| \\
&\leq & \frac 2{\sqrt \pi} \int^1_0 \big| \exp(-y^2 - 2uy \sqrt{-it} + itu^2)\big|du\\
&\leq & \frac 2{\sqrt \pi} \int^1_0  \exp( -y^2 -  \sqrt{2t} \ uy)du\\
&=& \sqrt{2/\pi} \ t^{-1/2}(1-\exp(- \sqrt{2t}\ y)) \exp(-y^2)/y\\
&<& t^{-1/2} \exp(-y^2)/y.
\end{eqnarray*}
From this a majorant is obtained for the above integral and theorem~3 is proved.
\end{proof}

\paragraph{\sc Remark. } In all the plotted examples the integrands in (\ref{eq:19}) exhibited a very smooth appearance. Since the integrand of the first integral in (\ref{eq:19}) tends to zero very rapidly the upper limit of integration can be lowered to moderate values for a numerical integration.

By Mir and Richards (1975) for $r \leq 1/\sqrt 2$ the elementary formula
\begin{equation}
\label{eq:21}
\Pr\{Q \leq r^2\} = \big( \sqrt n\ \pi^{(n-1)/2}/\Gamma((n+1)/2)\big)r^{n-1} - b(n)r^n, \ n \geq 3
\end{equation}
was derived where $b(3) = 2 \sqrt 6$ and $b(n)$ is given by an integral over the $(n-3)$--unit cube for $n >3$. The authors have computed $b(n)$ exactly for $n =3,4,5$ and additionally for $n = 6,7$ by numerical integration. Their values $b(n)$ have been checked by means of (\ref{eq:19}) and further values $b(n)$ can be computed now by univariate integration, e.g. $b(8) = 13.03951 \ldots$, $b(9) = 11.72490 \ldots$, $b(10) = 9.90095 \ldots$\ . 

The distribution of 
\[
(Q - \cE(Q))(\var(Q))^{-1/2} = (Q - (n-1)/12)\big((n-1)(2n+3)/(360 n)\big)^{-1/2}
\]
tends to a standard normal distribution for $n \to \infty$.

An example for the application of the distribution of $Q$ is the following statistical test. Let $Y$ be continuously distributed with a cdf $F$ and the median $m$. To test $F =F_0$, where the specified cdf $F_0$ has the same median, against
\[
F_0(y) < F(y), \ y<m \ \mbox{ and } \ F_0(y) > F(y), \ y>m ,
\]
we can use the test statistic $Q = \sum^n_{i=1} (U_i-\overline{U})^2$ with $U_i = F_0(Y_i)$, based on a random sample $Y_1,\ldots,Y_n$, and reject $F = F_0$ for large values of $Q$. In particular, with $F(y) = G((y-m)/\beta)$ and $F_0(y) = G((y-m)/\beta_0)$, $\beta > \beta_0 > 0$, where $G$ denotes a continuous (not normal) cdf with median$(G) = 0$, this is a test for the scale parameter $\beta$. It is supposed that the test works also if the medians of $F$ and $F_0$ differ only slightly, but the effect to the power is difficult to assess. The simpler test statistic $\sum^n_{i=1} (U_i - 1/2)^2$ with the exact cdf 
\[
2^{-(n-1)}\pi^{n/2-1} \sum^\infty_{k=1} \Im m\Big( \big( \erf(\sqrt{-ik\pi/n})/\sqrt{-ik\pi/n}\big)^n\Big) (1-\cos(4k\pi x/n))\big/ k,\ 0 \leq x < n/4
\]
can also be used.

%
\references


Abramowitz, M. and Stegun, I. (1968) \textit{Handbook of Mathematical Functions}, Dover Publications Inc., New York.

Bowman,  K.O. and Shenton, L.R. (1983) The distribution of the Standard Deviation and Skewness in Random Samples from a Gamma Density --- A New Look at a Craig--Pearson Study. Oak Ridge National Laboratory/CSD--109.

Box, G.E.P. (1953) Nonnormality and Tests on Variances, \textit{Biometrika} \textbf{40}, 318--335.

Craig, C.C. (1929) Sampling when the Parent Population is of Pearson's Type III,  \textit{Biometrika} \textbf{21}, 287--293.

Feller, W. (1971) \textit{An Introduction of Probability Theory and Its Applications}, Vol II, 2nd ed. John~Wiley \& Sons, New York. 

Mathai, A.M. and Provost, S.B. (1992) \textit{Quadratic Forms in Random Variables}, Marcel Dekker Inc., New York. 

Mir, M.A. and Richards, W.A. (1975) On a conjecture concerning the common content of an $n$--cube and a diagonal cylinder, \textit{Annals of the Institute of Statistical Mathematics} \textbf{27},  281--287.

Mudholkar, G.S. and Trivedi, M.C. (1981) A Gaussian Approximation to the Distribution of the Sample Variance for Nonnormal Populations, \textit{Journal of the American Statistical Association} \textbf{76}, 479--485.

Pearn, W.L. and Kotz, S. (2006) \textit{Encyclopedia and Handbook of Process Capability Indices: A Comprehensive Exposition of Quality Control Measures}. (Ser. on Quality, Reliability and Engineering Statistics, Vol. 12). World Scientific Publishing Company, Singapore. 

Pearn, W.L., Kotz, S. and Johnson, N.L. (1992) Distributional and Inferential Properties of Process Capability Indices, \textit{Journal of Quality Technology} \textbf{24}, 216--230.

Pearson, E.S. (1929) Note on Dr. Craig's paper, \textit{Biometrika} \textbf{21}, 294--302.

Roy, J. and Tiku, M.L. (1962) A Laguerre Series Approximation to the Sampling Distribution of the Variance, \textit{Sanky\=a}, Series A, \textbf{24}, 181--184.  

Royen, T. (2007a) Integral representations for convolutions of non--central multivariate gamma distributions, \mbox{arXiv:0704.0539} $\left[\textrm{math.ST}\right]$

Royen, T. (2007b) Exact distribution of the sample variance from a gamma parent distribution, \mbox{arXiv:0704.1415} $\left[\textrm{math.ST}\right]$

Tan, W.Y. and Wong, S.P. (1977) On the Roy--Tiku Approximation of Sample Variances from Non--normal Universes, \textit{Journal of the American Statistical Association} \textbf{71}, \mbox{875--880.}




\end{document}